\documentclass[reqno]{amsart}
\usepackage{amssymb,amsthm,amsfonts,amstext}
\usepackage{amsmath}
\usepackage{url}
\makeatletter \@addtoreset{equation}{section} \makeatother

\renewcommand\thefigure{\thesection.\@arabic\c@figure}
\renewcommand\thetable{\thesection.\@arabic\c@table}
\newtheorem{theorem}{Theorem}[section]
\newtheorem{lemma}[theorem]{Lemma}
\newtheorem{proposition}[theorem]{Proposition}

\newcommand{\mc}[1]{{\mathcal #1}}

\newcommand{\bb}[1]{{\mathbb #1}}

\renewcommand{\>}{\rangle}

 \DeclareMathOperator{\Var}{Var}

\title[Density fluctuations on percolation clusters]{Density fluctuations for a zero-range process on the percolation cluster}

\begin{document}

\author{Patr\'{\i}cia Gon\c{c}alves}
\address{Centro de Matem\'atica \\
Campus de Gualtar\\ 4710-057 Braga, Portugal} \email{patg@math.uminho.pt}
\author{ Milton Jara }
\address{FYMA\\ Universit\'e Catholique de Louvain\\Chemin du Cyclotron 2, B-1348\\ Louvain-la-Neuve, Belgium}
\email{milton.jara@uclouvain.be}

\date{\today}

\begin{abstract}
We prove that the density fluctuations for a zero-range process evolving on the supercritical percolation cluster
are given by a generalized Ornstein-Uhlenbeck process in the space of distributions $\mc S'(\bb R^d)$.
\end{abstract}


\renewcommand{\subjclassname}{\textup{2000} Mathematics Subject Classification}

\maketitle

\section{Introduction}

Consider the infinite cluster of a supercritical bond percolation model. On this random graph, we
define a zero-range process, which can be defined as a system of symmetric, simple random walks on which
the hopping time of particles at a given site depends only on the number of particles at that site. It has
 been recently proved that the scaling limit of a simple random walk on the percolation cluster is given by a
 Brownian motion with a diffusion coefficient that does not depend on the particular realization of the percolation
 cluster \cite{MP}, \cite{BB}. With this result in mind, it is natural to raise the question about the collective
 behavior of a system of random walks evolving on the percolation cluster, possibly with some interaction.
More precisely, we want to study the density fluctuations of this model. Scale the lattice by $1/n$ and give a mass
 $1/n$ to each particle. In this way we obtain a measure in $\bb R^d$, which we call the {\em empirical density} of particles.
  Under a diffusive time scaling, it has been proved for a simple exclusion process \cite{Fag} that the empirical measure converges
  in a proper sense, to a solution of the heat equation $\partial_t u = D\Delta u$, where $D$ is the diffusion coefficient of the underlying random walk.

This result can be interpreted as a law of large numbers for the density of particles, and therefore the question about the central limit
theorem arises. For the simple exclusion process, duality techniques allow to reduce the law of large numbers and the central limit theorem for
the density, to suitable problems for simple random walks. In order to treat a more general
 case, in these notes we discuss the density fluctuations for the zero-range process. Due to the inhomogeneity introduced by
 the randomness of the percolation cluster, the model turns out to be non-gradient. Following our previous work \cite{GJ}, a
 functional transformation of the density empirical allows to put the model back into the setup of gradient systems.
The main tool allowing to make this functional transformation is a form of the {\em compensated compactness lemma} of Tartar, which is obtained
in \cite{Fag} using two-scale convergence. Therefore, the proof of a central limit theorem for the density fluctuation is reduced to the
so-called {\em Boltzmann-Gibbs principle}, which roughly states that fluctuations of non-conserved quantities are faster than fluctuations of
conserved quantities. In the right scaling limit, only fluctuations of the density of particles survives, allowing to obtain a martingale
characterization of the scaling limit. In this way we prove that the limiting density field is given by a generalized Ornstein-Uhlenbeck
process, as predicted by the fluctuation-dissipation principle. Unfortunately, and only in dimension $d=2$, a regularity result for the
solutions of the heat equation on the percolation cluster prevent us to obtain such a result.

The main technical innovation in this paper is what we call the {\em connectivity lemma}. In \cite{GJ}, an ellipticity condition allows to
compare the relaxation properties of our system with the relaxation properties of a zero-range process in the absence of the random environment.
But the percolation cluster do not satisfy this ellipticity condition, since by construction the jump rate through a closed bond is 0. Moreover,
the percolation cluster is locally non connected, since for a typical box of fixed size, the part of the percolation cluster lying inside this
box is not connected. The connectivity lemma says that enlarging a bit the box we obtain a connected graph with good ergodic properties,
allowing to get bounds on the relaxation to equilibrium.

We point out that our results are still true in a more general context. We have chosen the zero-range process on the percolation cluster for simplicity and to capture the essential difficulties. Take for example a random barrier model on $\bb Z^d$ on which the conductances are not bounded below (see \cite{Fag} for a more detailed definition), and take an interacting particle system satisfying the gradient condition (before randomizing the lattice). Our results applies as well for these models, leading to the same results.

These notes are organized as follows. In Section \ref{s1} we introduce the model and we state the main results.
In Section \ref{s2} we introduce the corrected density field and we show how to obtain the main results from the compensated compactness lemma and the Boltzmann-Gibbs principle. In Section \ref{s3} we state and prove the connectivity lemma and we obtain the Boltzmann-Gibbs principle starting from it.

\section{Definitions and main results}
\label{s1}
\subsection{The supercritical percolation cluster}

Let $\mc E=\{e=\<xy\>;x,y \in \bb Z^d, |x-y|=1\}$ be the set of nearest-neighbor, non oriented bonds in $\bb Z^d$. Let $\omega = \{\omega(e); e \in \mc E\}$ be a sequence of i.i.d. random variables with $P(\omega(e)=1) = 1-P(\omega(e)=0)=p$. Let us fix a realization of the sequence $\omega$. Whenever $\omega(e)=1$, we say that the bond $e$ is {\em open}. Otherwise we say  that the bond $e$ is closed. We say that two sites $x, y$ in $\bb Z^d$ are {\em connected}, which we denote by $x \leftrightarrow y$, if there is a finite sequence of points $\{x_{0}=x,x_{1},\dots,x_{n}=y\}$ such that for any $i \in \{0,\dots,n-1\}$, $|x_{i+1}-x_i|=1$ and the bond $\<x_ix_{i+1}\>$ is open. For each $x \in \bb Z^d$, we define $\mc C(\omega,x) = \{y \in \bb Z^d; x \leftrightarrow y\}$. We call $\mc C(\omega,x)$ the cluster containing $x$. It is well known \cite{Grim} that the probability of $\mc C(\omega,x)$ being infinite is an increasing function of $p$, and that there exists $p_c \in (0,1)$ such that this probability is strictly positive for $p>p_c$ and equal to 0 for $p<p_c$. It is also well known that for $p > p_c$, with probability one, there exists a unique cluster  of infinite cardinality. We denote this cluster by $\mc C(\omega)$. We call $\mc C(\omega)$ the {\em supercritical percolation cluster}. From now on, we fix a number $p \in (p_c,1]$ and a sequence $\omega$ for which $\mc C(\omega)$ is well defined. We also define $\theta(p) = P(0 \in \mc C(\omega))$.

\subsection{The zero-range process in $\mc C(\omega)$}

For $x$, $y$ in $\bb Z^d$, we say that $x \sim y$ if $|x-y|=1$ and $\<xy\>$ is open. Let $g: \bb N_0 =\{0,1,...\} \to [0,\infty)$ be a function with $g(0)=0$. The zero-range process in $\mc C(\omega)$ with interaction rate $g(\cdot)$ is the Markov process $\xi_t$ in $\Omega =\bb N_0^{\mc C(\omega)}$ generated by the operator $L$ given by
\[
Lf(\xi) = \sum_{\substack{x, y \in \bb Z^d \\ x \sim y}} g(\xi(x)) \big[f(\xi^{x,y})-f(\xi)\big],
\]
where $f: \Omega \to \bb R$ is a {\em local} function, that is, it depends on $\xi(x)$ only for a finite number of sites $x \in \mc C(\omega)$, and $\xi^{x,y}$ is given by
\[
\xi^{x,y}(z) =
\begin{cases}
\xi(x)-1,&z=x\\
\xi(y)+1,&z=y\\
\xi(z),&z \neq x,y.\\
\end{cases}
\]

The dynamics of this process is easy to understand. At each time, a particle jumps from a site $x \in \mc C(\omega)$ to a neighbor $y \sim x$ with exponential rate $g(k)$, where $k$ is the number of particles at site $x$ at that time. When the initial number of particles is finite, the process $\xi_t$ is a continuous-time Markov chain in $\Omega$. For the construction of the process in the case of an infinite number of particles, we refer to \cite{And}.

In order to have a well defined family of ergodic, invariant measures for the process $\xi_t$, we assume that $g(n)>0$ for any $n>0$, and that
\begin{equation}
\label{A}
\quad \sup_{n} |g(n+1)-g(n)| < +\infty.
\end{equation}

Under these conditions, $g(\cdot)$ is bounded by a linear function: there exists $c_0$ such that $g(n) \leq c_0 n$ for any $n$. Let us define the uniform, product probability measures $\bar \nu_\varphi$ in $\Omega$ by the relation
\[
\bar \nu_\varphi\{ \xi(x) =k\} = \frac{1}{Z(\varphi)} \frac{\varphi^k}{g(k)!},
\]
where $g(k)!=g(1)\cdots g(k)$ for $k \geq 1$, $g(0)!=0$, $\varphi \in [0,\infty)$ and $Z(\varphi)$ is the normalization constant. Due to (\ref{A}), these measures are well defined up to some critical value $\varphi_c$. Notice that the parameter $\varphi$ has an interpretation as the strenght of the interaction, since $\int g(\xi(x))\bar \nu_\varphi(d\xi)=\varphi$. It is not hard to prove that the measures $\{\bar \nu_\varphi; \varphi \in [0,\varphi_c)\}$ are ergodic and invariant under the evolution of $\xi_t$. Since we are interested in density fluctuations, we will reparametrize this family of measures by its density of particles.
For that purpose, define $\rho(\varphi)= \int \xi(x) \bar \nu_\varphi(d\xi)$. It is not hard to see that $\varphi \mapsto \rho(\varphi)$ is a diffeomorphism from $[0,\varphi_c)$ to $[0,\infty)$. In particular, the inverse function $\rho \mapsto \varphi(\rho)$ is well defined. We define $\nu_\rho = \bar \nu_{\varphi(\rho)} $, measure for which now $\int \xi(x) \nu_\rho(d\xi)=\rho$.

\subsection{The density fluctuations}

Let $\rho \in (0,\infty)$ be fixed and consider the process $\xi_t^n= \xi_{tn^2}$ starting from the measure $\nu_\rho$. We denote by $\bb P_n$ the distribution of the process $\xi_\cdot^n$ in the Skorohod space $\mc D([0,\infty),\Omega)$ of c\`adl\`ag trajectories in $\Omega$, and by $\bb E_n$ the expectation with respect to $\bb P_n$.

Denote by $\mc C_c(\bb R^d)$ the set of continuous functions $G: \bb R^d \to \bb R$ with compact support. By the ergodic theorem,
 it is not hard to see that
\[
\lim_{n \to \infty} \frac{1}{n^d} \sum_{x \in \mc C(\omega)} G(x/n) \xi(x) = \theta(p) \rho \int G(x) dx,
\]
almost surely with respect to the probability measure $P \otimes \nu_\rho$. In the previous result, the ergodic theorem is invoked twice: first
to state that the density of points belonging to $\mc C(\omega)$ when properly rescaled is equal to $\theta(p)$, and then to state that the
density of particles is equal to $\rho$. Since the measure $\nu_\rho$ is invariant under the dynamics of $\xi$, the same result is also valid if
we replace $\xi$ by $\xi_t^n$ in the previous expression. We are interested on a version of the central limit theorem for this quantity. Let us
define, for each test function $G$, the {\em density fluctuation field} $\mc Y_t^n(G)$ by
\[
\mc Y_t^n(G) = \frac{1}{n^{d/2}} \sum_{x \in \mc C(\omega)} G(x/n) \big(\xi_t^n(x)-\rho\big).
\]

For topological reasons, it will be convenient to restrict the previous definition to functions $G \in \mc S(\bb R^d)$, the Schwartz space of test functions, although $\mc Y_t^n(G)$ makes sense for more general test functions. The process $\mc Y_t^n$ defined in this way corresponds to a process on the Skorohod space $\mc D([0,\infty),\mc S'(\bb R^d))$, where $\mc S'(\bb R^d)$ is the space of distributions on $\bb R^d$. Now we are ready to state our main result:

\begin{theorem}
\label{t1}
Fix a particle density $\rho>0$. In dimension $d>2$, for almost all $\omega$, the sequence of processes $\{\mc Y_t^n\}_n$ converges in the sense of finite-dimensional distributions to a generalized Ornstein-Uhlenbeck process $\mc Y_t$ of mean zero and characteristics $\varphi'(\rho) D \Delta$, $\sqrt{D\varphi(\rho)}\nabla$, where $D$ is the limiting variance of a symmetric random walk in $\mc C(\omega)$.
\end{theorem}

Notice that this theorem implies in particular a central limit theorem for the sequence $\{\mc Y_t^n(G)\}$ for any $G \in \mc S(\bb R^d)$.

\section{The corrected fluctuation field}
\label{s2}

\subsection{The corrected fluctuation field}

The proof of Theorem \ref{t1} is based on Holley-Stroock's characterization of generalized Ornstein-Uhlenbeck processes \cite{HS}.
For each $t \geq 0$, let $\mc F_t$ be the $\sigma$-algebra on $\mc D([0,\infty),\mc S'(\bb R^d))$ generated by the projections $\{\mc Y_s(H); s\leq t,  H \in \mc S(\bb R^d)\}$. The process $\mc Y_t$ admits the following characterization:

\begin{proposition}
\label{p1}
There exists a unique process $\mc Y_t$ in $\mc C([0,\infty),\mc S'(\bb R^d))$ such that:
\begin{itemize}
\item[i)] For every function $G \in \mc S(\bb R^d)$,
\[
M_t(G)= \mc Y_t(G) -\mc Y_0(G) - \varphi'(\rho) D \int_0^t \mc Y_s(\Delta G)ds
\]
and
\[
\big(M_t(G)\big)^2 - \theta(p)\varphi(\rho)Dt \int_{\bb R^d} \big(\nabla G(x)\big)^2 dx
\]
are $\mc F_t$-martingales.

\item[ii)] $\mc Y_0$ is a Gaussian field of mean zero and covariance given by
\[
E\big[ \mc Y_0(G) \mc Y_0(H)\big] = \theta(p) \chi(\rho) \int_{\bb R^d} G(x) H(x) dx,
\]
where $\chi(\rho) = \Var(\eta(0); \nu_\rho)$ and $G,H \in \mc S(\bb R^d)$. The process $\mc Y_t$ is called the generalized Ornstein-Uhlenbeck process of mean-zero and characteristics $\varphi'(\rho) D \Delta$ and $\sqrt{\varphi(\rho) D} \nabla$.
\end{itemize}
\end{proposition}

Now the idea behind the proof of Theorem \ref{t1} is simple. We will prove that the sequence of processes $\{ \mc Y_t^n\}_n$ is tight and that every limit point satisfies the martingale problem stated in Proposition \ref{p1}.  With these two elements in hand, we will be able to conclude Theorem \ref{t1}.
By Dynkin's formula, for each $G \in \mc S(\bb R^d)$,
\[
M_t^n(G) = \mc Y_t^n(G) -\mc Y_0^n(G) - \int_0^t n^2 L \mc Y_s^n(G) ds
\]
is a martingale with respect to $\mc F_t^n = \sigma\{\xi_s^n;s \leq t\}$. The quadratic variation of $M_t^n(G)$ is given by
\[
\int_0^t n^2\big\{ L \mc Y_s^n(G)^2 - 2\mc Y_s^n(G) L\mc Y_s^n(G)\big\}ds.
\]

Notice that the second step (that is, that the limit points satisfy the martingale problem) requires to replace, in some sense, $n^2 L \mc Y_s^n(G)$ by $\varphi'(\rho) D \mc Y_s^n(\Delta G)$, and to replace $L \mc Y_s^n(G)^2 - 2\mc Y_s^n(G) L\mc Y_s^n(G)$ by $\theta(p) \varphi(\rho) D\int_{\bb R^d} \big(\nabla G(x)\big)^2 dx$. Let us take a more careful look at these two terms. We start with the second one. Simple computations show that
\[
L \mc Y_s^n(G)^2 - 2\mc Y_s^n(G) L\mc Y_s^n(G) = \frac{1}{n^d} \sum_{\substack{x \in \mc C(\omega)\\y:x \sim y}}  g\big(\xi_s^n(x)\big)n^2 \big(G(y/n)-G(x/n)\big)^2.
\]

By the ergodic theorem, when integrated in time this quantity should converge to $\theta(p) \varphi(\rho) t \int \big(\nabla G(x)\big)^2 dx$. Notice that we have missed the factor $D$ in the previous computation. We will return to that point later. Now let us take a look at the first term:
\[
n^2 L \mc Y_s^n(G) = \frac{1}{n^{d/2}} \sum_{x \in \mc C(\omega)} \big(g(\xi_s^n(x))-\varphi(\rho)\big) \mc L_n G(x/n),
\]
where $g(\cdot)$ is the interaction rate and $\mc L_n $ is the generator of the associated symmetric random walk in $\mc C(\omega)$. For a function $F:\bb R^d \to \bb R$ and $x \in \mc C(\omega)$, $\mc L_n F(x/n)$ is defined by
\[
\mc L_n F(x/n) = n^2\sum_{\substack{y \in \mc C(\omega)\\ y \sim x}} F(y/n) -F(x/n).
\]

The so-called {\em Boltzmann-Gibbs principle} will allow us to replace, when integrated in time, the expression $g(\xi_s^n(x)) -\varphi(\rho)$ by $\varphi'(\rho)\big(\xi_s^n(x)-\rho\big)$ in the sum above, allowing us to write $n^2 L \mc Y_t^n(G)$ as $\varphi'(\rho) \mc Y_t(\mc L_n G)$ plus a term that is negligible as $n \to \infty$.
Now we can see what the problem is. Assume that $\{\mc Y_t^n\}_n$ is tight. Take a limit point $\mc Y_t^\infty$ of $\{\mc Y_t^n\}_n$. We want to say that $\mc Y_t^n(\mc L_n G)$ converges to $\mc Y_t^\infty(D\Delta G)$ along the corresponding subsequence. But $\mc Y_t^n$ converges to $\mc Y_t^\infty$ only on a weak sense. Therefore, we should need strong convergence of $\mc L_n G$, which is easily checked not to hold.

The way to overcome this problem is to use the {\em compensated compactness lemma}. The idea is to choose, for each $n$, a test function $G_n$ in such a way that $\mc L_n G_n$ converges strongly to $D\Delta G$ as $n \to \infty$. We therefore define the {\em corrected fluctuation field} as in \cite{GJ}. Fix some $\lambda >0$. For each $G \in \mc S(\bb R^d)$, define $G_n: \mc C(\omega) \to \bb R^d$ as the solution of the {\em resolvent} equation
\begin{equation}
\label{ec1}
\lambda G_n(x) - \mc L_n G_n(x) = \lambda G(x/n) - D \Delta G(x/n).
\end{equation}

Then the corrected fluctuation field $\mc Y_t^{n,\lambda}$ is given by
\[
\mc Y_t^{n,\lambda}(G) = \frac{1}{n^{d/2}} \sum_{x \in \mc C(\omega)} \big(\xi_t^n(x) -\rho\big) G_n(x).
\]

Notice that $\mc L_n G_n = \lambda(G_n - G) + D \Delta G$. In particular, strong convergence of $\mc L_n G_n$ to $D \Delta G$ follows from strong convergence of $G_n$ to $G$. The following proposition tells us that this is, indeed, the case.

\begin{proposition}[Faggionato \cite{Fag}]
\label{p2}
There is a set of $P$-total probability such that for any $G \in \mc S(\bb R^d)$,
the sequence $\{G_n\}_n$ converges to $G_n$ in the following strong sense:
\[
\lim_{n \to \infty} \frac{1}{n^d} \sum_{x \in \mc C(\omega)} \big|G_n(x) -G(x/n)\big|^2 =0.
\]
\end{proposition}

In particular, since the invariant measure $\nu_\rho$ is of product form, we conclude that $\mc Y_t^{n,\lambda}(G_n) -\mc Y_t^n(G)$ vanishes in $\mc L^2(\bb P_n)$ as $n \to \infty$. We will see that the standard scheme, tightness plus uniquenes of limit points via martingale characterization, can be accomplished for $\mc Y_t^{n,\lambda}$.

\subsection{The martingale problem for the corrected fluctuation field}

When considering the corrected fluctuation field $\mc Y_t^{n,\lambda}$ instead of $\mc Y_t^n$, the martingale representation gives
\[
M_t^{n,\lambda}(G) = \mc Y_t^{n,\lambda}(G) - \mc Y_0^{n,\lambda}(G)
    - \int_0^t \Theta_s^n(\lambda(G_n-G) +D\Delta G) ds.
\]

In this formula, we have defined
\[
\Theta_t^n(F) = \frac{1}{n^{d/2}} \sum_{x \in \mc C(\omega)}  \big(g(\xi_t^n(x))-\varphi(\rho)\big) F(x/n).
\]

The quadratic variation of $M_t^{n,\lambda}$ is given by
\[
\<M_t^{n,\lambda}\> = \int_0^t \frac{1}{n^d} \sum_{\substack{x \in \mc C(\omega)\\y: x \sim y}} n^2 g\big(\xi_s^n(x)\big)\big(G_n(y)-G_n(x)\big)^2ds.
\]

Multiplying the resolvent equation (\ref{ec1}) by $G_n(x)$ and summing over $x \in \mc C(\omega)$ we see that
\begin{equation}
\label{ec2}
\lim_{n \to \infty}
\frac{1}{n^d} \sum_{\substack{x \in \mc C(\omega)\\y: x \sim y}} n^2
    \big(G_n(y)-G_n(x)\big)^2 = \theta(p)D \int\big(\nabla G(x)\big)^2 dx.
\end{equation}

Let us denote by $Q_n^\lambda$ the distribution of the process $\mc Y_\cdot^{n,\lambda}$ in $\mc D([0,\infty), \mc S'(\bb R^d))$.
Following \cite{GJ}, it is not hard to prove that the sequence $\{Q_n^\lambda\}$ is tight. In the same way, for $d \geq 3$, we can prove that
\[
\lim_{n \to \infty} \<M_t^{n,\lambda}\> = \varphi(\rho) \theta(p) D \int \big(\nabla G(x)\big)^2dx.
\]

The arguments in \cite{GJ}, and only in dimension $d=2$, require some regularity of $G_n$ that is missing in our situation. This is the only point on this article where we require $d \geq 3$, otherwise the proofs do not depend on the dimension. Now we will state in a more precise way the Boltzmann-Gibbs principle.

\begin{theorem}[Boltzmann-Gibbs principle]
\label{t2} For any $t>0$ and any $F \in \mc C_c(\bb R^d)$,
\[
\lim_{n \to \infty} \bb E_n\Big[ \Big(\int_0^t \big\{\Theta_s^n(F) - \varphi'(\rho) \mc Y_s^n(F)\big\} ds\Big)^2\Big] =0.
\]
\end{theorem}

As in \cite{GJ}, from the tightness of $\{Q_n^\lambda\}$, Theorem \ref{t2}  an the expression (\ref{ec2}) allow us conclude the following result:

\begin{theorem}
\label{t3}
Fix a particle density $\rho>0$. In dimension $d>2$, for almost any $\omega$, the sequence of processes $\{\mc Y_t^{n,\lambda}\}_n$ converges in distribution with respect to the Skorohod topology in $\mc D([0,\infty),\mc S'(\bb R^d))$ to a generalized Ornstein-Uhlenbeck process $\mc Y_t$ of mean zero and characteristics $\varphi'(\rho)D \Delta$, $\sqrt{D\varphi(\rho)}\nabla$, where $D$ is the limiting variance of a symmetric random walk in $\mc C(\omega)$.
\end{theorem}

Theorem \ref{t1} is an immediate consequence of this Theorem and Proposition \ref{p2}.

\section{The Boltzmann-Gibbs principle}
\label{s3}

In this section we prove Theorem \ref{t2}. The main point that makes the proof of this theorem different from the proof in \cite{GJ} is the lack
of ellipticity: by construction, the jump rate between two neighboring sites $x$, $y$ is equal to 0 when the bond $\<xy\>$ is not open. In
\cite{GJ}, due to the ellipticity condition we could compare the generator of the process with the generator of a zero-range process without the
random environment, but with a slower jump rate. Since the Boltzmann-Gibbs principle holds for the latter process, it should hold for the
former, since the dynamics is faster. In our case, we do not have a proper slower process to compare with.

Let us explain better the intuition behind Theorem \ref{t2}. The idea is that non-conserved quantities fluctuate faster than conserved ones.
Since the only conserved quantity for the zero-range process is the number of particles, it is reasonable that at the right scale, the only part
of the fluctuation field $\Theta_t^n$ that is seen at a macroscopic level is its projection over the conservative field $\mc Y_t^n$. We can
think that non-conserved quantities equilibrate locally, while conserved quantities need to be transported in order to equilibrate. In
particular, we will see that the form of the graph is not really important in Theorem \ref{t2}. What is really important is the connectivity of
the graph: if a graph has more than one connected component, then there is more that one conserved quantity: the number of particles on each
connected component.

For $n \in \bb N_0$, define $\mc V(n) = g(n) -\varphi(\rho) -\varphi'(\rho)(n-\rho)$. The strategy of the proof of Theorem \ref{t2} introduced
by Chang \cite{Cha} is the following. Fix a positive integer $k$. Divide the support of the test function $F$ into small boxes of size $k/n$.
Since $F$ is continuous, we can average the function $\mc V(\xi(x))$ over the corresponding boxes on the lattice. Then we use some sort of
ergodic theorem to reduce the sum over many blocks integrated in time, into a sum over a single block. This last problem is a static one. A new
ergodicity argument, now with respect to the invariant measure $\nu_\rho$ will allow us to conclude.

An important property of the lattice is that it is {\em locally connected}, in the sense that the restriction of the lattice to any box is still a connected graph. This is no longer true for the percolation cluster: the restriction of the percolation cluster to a box is, in general, not connected. Although most of it belongs to a big single component, the rest is spread over many small connected components. Of course, all these connected components are connected by paths that pass outside the initial box. The point is that these paths can be chosen is such a way that they do not go too far from the original box.
We will develop these ideas in the next paragraphs.

\subsection{A connectivity lemma}

In this section we will state a result that we call the {\em connectivity lemma} and we will prove Theorem \ref{t2} starting from it. The proof of this connectivity lemma is postponed to the next section.

Let $k$, $l$ be fixed positive integers. We will send $k$ and $l$ to infinity after $n$. We introduce two intermediate scales in our problem as follows. For simplicity we assume $n/k \in \bb N$, and pasting a sufficient number of cubes, we can assume that the support of $F$ is contained in the cube $(\delta,1-\delta)^d$ for some $\delta>0$. We can split the cube $\Lambda_n=\{1,2,\dots,n\}^d$ into $(n/k)^d$ non-overlapping cubes of side $k$. Let $\{\bar B_j^0,j=1,\dots,(n/k)^d\}$ be an enumeration of those cubes. Define $\bar B_j$ as the box of side $(2l+1)k$, centered at $\bar B_j^0$.
In particular, the cubes $\bar B_j$ are the union of $(2l+1)^d$ cubes in $\{\bar B_i^0\}_i$, except for the ones near the border of $\Lambda_n$.
Now define
\[
B_j^0 = \bar B_j^0 \cap \mc C(\omega),
\]
\[
B_j = \{x \in \bar B_j; x \leftrightarrow B_j^0\}.
\]

In other words, $B_j^0$ is the portion of the cluster $\mc C(\omega)$ inside the box $\bar B_j^0$, and $B_j$ is the portion of the cluster $\mc C(\omega)$ inside $\bar B_j$ and connected to the box $\bar B_j^0$. We say that a cube $\bar B_j^0$ is {\em good} if $B_j$ is connected. In other words, $\bar B_j^0$ is good if the connected components of $B_j^0$ are connected between them by paths that lie entirely in $\bar B_j$.
We will denote by $\mc B_n$ the union of good cubes. The numbers $k$ and $l$ will be fixed most of the time.  Therefore, in order to keep notation simple we do not make explicit the dependence of $\mc B_n$ in $k$ and $l$.
 A cube that is not in $\mc B_n$ will be called {\em bad}. We will call $B_j$ indistinctely the set of points already defined, and the subgraph of $\mc C(\omega)$ corresponding to these points.

Define $\mc C_n(\omega) = \mc C(\omega) \cap \Lambda_n$.
Theorem \ref{t2} is a consequence of the following lemma:

\begin{lemma}[Connectivity lemma]
For each $\epsilon>0$ there exists $l>0$ such that
\begin{align}
\label{ec2.1}
    i)
    &\lim_{k \to \infty} \limsup_{n \to \infty} \bb E_n \Big[ \Big(\int_0^t \frac{1}{n^{d/2}} \sum_{x \in \mc B_n} \mc V(\xi_s^n(x)) F(x/n)ds\Big)^2 \Big] =0\\
    \label{ec2.2}
    ii)
    &\limsup_{k \to \infty} \limsup_{n \to \infty} \frac{|\mc C_n(\omega) \setminus \mc B_n|}{n^d} \leq \epsilon,
\end{align}
where $|\mc C_n(\omega) \setminus \mc B_n|$ denotes the cardinality of the set $\mc C_n(\omega) \setminus \mc B_n$.
\label{l1}
\end{lemma}

The first part of the lemma says, roughly speaking, that Theorem \ref{t2} is true if the considered graph is locally connected.
The second part says that $\mc C(\omega)$ is not locally connected only on a small portion of the lattice.

Now let us prove Theorem \ref{t2} assuming Lemma \ref{l1}. The expectation appearing in the statement of Theorem \ref{t2} is bounded by
\begin{multline*}
2\bb E_n\Big[ \Big(\int_0^t \frac{1}{n^{d/2}} \sum_{x \in \mc B_n} \mc V(\xi_s^n(x)) F(x/n) ds\Big)^2 \Big]\\
    +2\bb E_n\Big[ \Big(\int_0^t \frac{1}{n^{d/2}} \sum_{x \in \mc C_n(\omega) \setminus \mc B_n} \mc V(\xi_s^n(x)) F(x/n) ds\Big)^2 \Big].
\end{multline*}

The first expectation goes to 0 as $n \to \infty$ and then $k \to \infty$ by Lemma \ref{l1}, part $i)$. By Schwartz inequality, the second expectation is bounded by
\[
\frac{t^2}{n^d} \sum_{x \in \mc C_n(\omega) \setminus \mc B_n} F(x/n)^2 E_\rho[\mc V(\xi(x))^2],
\]
which turns out to be bounded by $C(t,F,\rho)|\mc C_n(\omega) \setminus \mc B_n|/n^d$. This last expression vanishes as $n \to \infty$ and then $k \to \infty$ by Lemma \ref{l1}, part $ii)$, which proves Theorem \ref{t2}.

\subsection{Proof of the connectivity lemma, part $i)$}
In order to simplify the notation, we also denote by $\mc B_n$ the set of indices $j$ in $\{1,\dots,(n/k)^d\}$ for which $\bar B_j^0$ is a good cube.
Let us start with part $i)$. We will manipulate the term
\begin{equation}
\label{ec3}
\bb E_n \Big[ \Big(\int_0^t \frac{1}{n^{d/2}} \sum_{x \in \mc B_n} \mc V(\xi_s^n(x)) F(x/n)ds\Big)^2 \Big]
\end{equation}
until we arrive to an expression that does not depend on $t$. Taking the positive and negative parts of $F$, we can assume, without loss of
generality, that $F$ is non-negative. For each $j$, take a point $y_j$ in $\bar B_j^0$. Since the function $F$ is uniformly continuous, we can
rewrite the integrand in (\ref{ec3}) as
\begin{equation}
\label{ec4}
\frac{1}{n^{d/2}} \sum_{j \in \mc B_n} F(y_j/n) \sum_{x \in \bar B_j^0} \mc V(\xi_s^n(x))
\end{equation}
plus a rest that vanishes in $\mc L^2(\nu_\rho)$ as $n \to \infty$ and then $k \to \infty$. Notice that now we have introduced an averaging of the function $\mc V(\xi(x))$ over boxes of side $k$.

Now we explain the point where we make use of the time average in (\ref{ec3}). For each $j$, denote by $L_{B_j}$ the restriction of the
generator $L$ to the set $B_j$. Observe that the zero-range process restricted to $\bb N_0^{B_j}$ is ergodic on the set of configurations with a
fixed number of particles, exactly due to the fact that $B_j$ is a good cube. For any two functions $f, h: \Omega \to \bb R$, denote by
$\<f,h\>_\rho$ the inner product with respect to $\nu_\rho$. We define by
\[
||f||_1^2 = \<f,-Lf\>_\rho,
\]
\[
||f||_{-1}^2 = \sup_h\{ 2\<f,h\>_\rho -||h||_1^2\}
\]
the Sobolev norms associated to $L$. Here the supremum is taken over functions $h \in \mc L^2(\nu_\rho)$.
Take an arbitrary (by now) family of functions $\{f_j;j \in \mc B_n\}$ with $f_j :\bb N_0^{B_j} \to \bb R$.
We have the following {\em Sobolev inequality} (Prop. A1.6.1 of \cite{KL}):
\[
\bb E_n\Big[\Big(\int_0^t \frac{1}{n^{d/2}} \sum_{j \in \mc B_n} F(y_j/n) L_{B_j} f_j(\xi_s^n)ds\Big)^2\Big]
    \leq \frac{20 t}{n^2} \big|\big|\frac{1}{n^{d/2}} \sum_{j \in \mc B_n} F(y_j/n) L_{B_j} f_j\big|\big|_{-1}^2.
\]

Let us call $\mc W_n$ the term inside the norm.  We will bound $||\mc W_n||_{-1}^2$ using the variational formula of $||\cdot||_{-1}^2$ introduced above. Then, we need to estimate $\<h,\mc W_n\>_\rho$ for $h$ given. By the weighted Schwartz inequality,
\[
\int h L_{B_j} f d\nu_\rho \leq \frac{1}{2\gamma_j} \<f,-L_{B_j} f\>_\rho +\frac{\gamma_j}{2} \<h,-L_{B_j} h\>_\rho.
\]

Choose $\gamma_i= n^{2+d/2}/F(y_j/n)$ and observe that $\sum_j \<h,-L_{B_j}h\>_\rho \leq \<h,-Lh\>_\rho$. Plug
 this estimate into the variational formula for $||\mc W_n||_{-1}^2$ to discover that
\[
\bb E_n\Big[\Big(\int_0^t \frac{1}{n^{d/2}} \sum_{j \in \mc B_n} F(y_j/n) L_{B_j} f_j(\xi_s^n)ds\Big)^2\Big]
    \leq \frac{40 t}{n^{d+2}}\sum_{j \in \mc B_n} F(y_j/n)^2 \<f_j,-L_{B_j}f\>_\rho.
\]

Notice that, for $f_j$ fixed, the right-hand side of this inequality is of order $1/n^2$. Notice also that we can not take the same function
$f_j$ for different $j$, since the clusters $B_j$ are different for each $j$. However, $k$ and $l$ are fixed. In particular, the possibilities
for the cluster $B_j$ are many, but finite. Let $\Xi$ be the set of connected graphs contained in a box of size $(2l+1)k$. Let us denote by
$\eta$ a generic element (that is, a graph) in $\Xi$. Take now a family of functions $\{f^\eta; \eta \in \Xi\}$ and define $f_j = f^\eta$ if
$B_j = \eta$. With this notation, we can bound the right-hand side of the previous inequality by
\[
\frac{40 t||F||_\infty^2}{n^{d+2}}\sum_{\eta \in \Xi} q^n(\eta) \<f^\eta,-L_{\eta}f^\eta\>_\rho,
\]
where $q^n(\eta)$ is the number of boxes in $\mc B_n$ with graph $\eta$, and $L^\eta$ is the generator $L$ restricted to the graph $\eta$. Since $q^n(\eta)/n^d$ is bounded, the functions $\{f^\eta\}_\eta$ are fixed and also $\Xi$ is fixed, this last quantity goes to 0 as $n \to \infty$.
In particular, this means that we can discount a sum of the form $n^{-d/2} \sum_j F(y_j/n) L_{B_j} f_j$ in (\ref{ec4}). We can also take the infimum over the families $\{f^\eta\}_\eta$, although only after $n \to \infty$. In other words, we have reduced (\ref{ec3}) to the verification of
\[
\lim_{k \to \infty} \sup_{\{f^\eta\}_\eta} \limsup_{n \to \infty}
    \bb E_n \Big[ \Big(\int_0^t \frac{1}{n^{d/2}} \sum_{j \in \mc B_n} F(y_j/n) \sum_{x \in \bar B_j^0} \big\{\mc V(\xi_s^n(x)) - f_j(\xi_s^n)\big\}
ds\Big)^2 \Big].
\]

The good point is that in this sum the time integral has already played its part, allowing the introduction of the functions $f_j$.
Bound the expectation using Schwartz inequality twice, once to get rid of the time integral, and once more to make use of the product form of the measure $\nu_\rho$.
Notice that the boxes $B_j$ at a distance greater than $(2l+1)k$ are independent. Then, the expectation in the previous limit is bounded by
\[
\frac{t^2(2l+1)^d||F||_\infty^2}{n^d}\sum_{\eta \in \Xi} q^n(\eta) E_\rho\Big[\Big(\sum_{x \in \eta} \mc V(\xi(x)) - L^\eta f^\eta\Big)^2\Big].
\]

Remark that the sum above depends on $n$ only through the number $q^n(\eta)/n^d$. By the ergodic theorem, $q^n(\eta)/(n/k)^d$ converges, as $n \to \infty$, to $q(\eta)$, which is the probability of a given subgraph $B_j$ of $\mc C(\omega)$ be (good and) equal to $\eta$. Therefore, we are left to prove that
\begin{equation}
\label{ec5}
\lim_{k \to \infty} \sup_{\{f^\eta\}_\eta}
\frac{t^2(2l+1)^d||F||_\infty^2}{k^d}\sum_{\eta \in \Xi} q(\eta) E_\rho\Big[\Big(\sum_{x \in \eta} \mc V(\xi(x)) - L^\eta f^\eta\Big)^2\Big]=0.
\end{equation}

Here we remark that the whole construction we have done has as purpose to get a {\em connected} graph in the sum above. For each graph $\eta$ fixed, and due to the ergodicity of the zero-range process on sets with fixed number of particles,
\[
\inf_{f^\eta} E_\rho \Big[\Big(\sum_{x \in \eta} \mc V(\xi(x)) - L^\eta f^\eta\Big)^2\Big]
    = E_\rho \Big[E_\rho\Big(\sum_{x \in \eta} \mc V(\xi(x))\Big| \xi^\eta\Big)\Big],
\]
where $\xi^\eta$ denotes the number of particles of the configuration $\xi$ on the graph $\eta$. Define, for $\rho' \geq 0$, $\psi(\rho') = E_{\rho}[\mc V(\xi(x))]$. Notice that $\psi(\rho)= \psi'(\rho)=0$. This fact combined with the equivalence of ensembles allow us to show that
\[
\sup E_\rho \Big[E_\rho\Big(\sum_{x \in \eta} \mc V(\xi(x))\Big| \xi^\eta\Big)\Big] <+\infty,
\]
where the supremum is over $k>0$ {\em and} $\eta \in \Xi$. Therefore, we have proved that the expression in (\ref{ec5}) is bounded by
\[
\frac{C(F,g) t^2 (2l+1)^d}{k^d} \sum_{\eta \in \Xi} q(\eta),
\]
which goes to 0 as $k \to \infty$, since the sum in $\eta$ is bounded by 1 (remember that $q(\eta)$ is a probability).

\subsection{Proof of the connectivity lemma, part $ii)$}
We will use a result of \cite{AP},  which roughly states that the percolation cluster has good connectivity properties.

For any two points $x, y \in \mc C(\omega)$, we define the distance $D(x,y)$ as the length of the minimal path connecting $x$ and $y$:
\[
D(x,y)  = \inf\{n; \exists x=x_0,\dots,x_n=y \text{ with } \<x_{i-1}x_i\> \text{ open}\}.
\]

We have the following result:

\begin{proposition}[Antal-Pisztora \cite{AP}]
There exists a constant $\gamma= \gamma(p,d) \in [1,\infty)$ such that
\[
\limsup_{|z| \to \infty} \frac{1}{|z|} \log P\big(0\leftrightarrow z , D(0,z) > \gamma|z|\big) <0.
\]
\label{p3}
\end{proposition}

The previous Proposition is telling us that, with high probability, two points that belong to $\mc C(\omega)$, both inside a box of size $k$, are connected through a path contained on the box of size $(\gamma +\delta)k$, centered on the original box. This means that for $l>\gamma$, (\ref{ec2.2}) should hold. Now we will make this point more precise.

For simplicity, assume $k$ odd. Fix $\epsilon>0$.
Let $\Gamma_k^0$ be the box of size $k$, centered at the origin.
Define $\Gamma_k= \Gamma_k^0 \cap \mc C(\omega)$, the intersection of $\Gamma_k^0$ with the infinite cluster. By the ergodic theorem, $|\Gamma_k|/k^d \to \theta(p)$ $a.s.$ and also in probability.
Fix $k_0$ such that $\Gamma_{k_0} \neq \emptyset$ with probability bigger that $1-\epsilon/2$.
By Proposition \ref{p3}, there is $\delta>0$ such that
\begin{equation}
\label{ec6}
P(0 \leftrightarrow z, D(0,z) >\gamma|z|) \leq \exp\{-\delta|z|\}
\end{equation}
for $|z|$ is big enough. Taking a bigger $k_0$ if necessary, we can assume that this inequality holds for any $z$ such that $|z|>k_0$. Fix $x \in \Gamma_{k_0}^0$ and consider for a moment the measure $P(\cdot|x \in \mc C(\omega))$, well defined since $P(x \in \mc C(\omega))=\theta(p)>0$. Fix $l>0$. In analogy with our previous definition, we say that the box $\Gamma_k$ is good if each connected component of $\Gamma_k$ is connected to $x$ through a path contained in $\Gamma_{(2l+1)k}^0$.
This is the case if and only if each point of $\mc C(\omega) \cap \partial \Gamma_k^0$ is connected to $x$ by a path contained in $\Gamma_{(2l+1)k}^0$.
When we do not fix a point $x \in \Gamma_{k_0}^0$, we say that $\Gamma_k$ is good if any two connected components of $\Gamma_k$ are connected by a path contained in $\Gamma_{(2l+1)k}^0$.
For $k>k_0$, each one of the points in $\mc C(\omega) \cap \partial \Gamma_k^0$ are at distance  from $x$ at least $|k-k_0|$ and at most $|k+k_0|$. By (\ref{ec6}), the probability (with respect to $P(\cdot|x \in \mc C(\omega))$) of all these points being connected to $x$ by a path contained in $\Gamma_{(2l+1)k}^0$ is bigger than
\begin{equation}
\label{ec7}
1- 2d k^{d-1} \exp\{-\delta |k-k_0|\}
\end{equation}
as soon as $l \geq \gamma +(1+\gamma)k_0/k$. This is the case if, for example, $l \geq 3\gamma$. The quantity in $(\ref{ec7})$ goes to 0 as $k \to \infty$. Since this estimate is independent of the choice of site $x \in \Gamma_{k_0}^0$, we conclude that
\[
\limsup_{k \to \infty} P(\Gamma_k \text{ is not good }) \leq \epsilon/2.
\]

Notice that the definition of {\em goodness} of a given box is translation invariant. In particular, with the $B_j$'s notation,
\[
\limsup_{k \to \infty} P(B_j \text{ is not good }) \leq \epsilon/2.
\]
By the ergodic theorem,
\[
\lim_{n \to \infty} \frac{|\mc C_n(\omega) \setminus \mc B_n|}{n^d} = P(B_j \text{ is not good})
\]
in probability and $a.s.$ We conclude that for $l > 3\gamma$,
\[
\limsup_{k \to \infty} \limsup_{n \to \infty} \frac{|\mc C_n(\omega) \setminus \mc B_n|}{n^d} \leq \epsilon,
\]
which ends the proof of (\ref{ec2.2}).

\section*{Acknowledgements}
M.J. was supported by the Belgian Interuniversity Attraction Poles Program P6/02,
through the network NOSY (Nonlinear systems, stochastic processes and statistical mechanics).
 M.J. would like to thanks the hospitality of Universidade do Minho, where part of this work was done.

P.G. was supported by F.C.T. Portugal with the grant SFRH/BPD/39991/2007. P.G. wants to thank the hospitality of the Universit\'e Catholique de
Louvain-la-Neuve and IMPA, where part of this work was done.

 The authors would like to thank A. Faggionato for giving access to ref. \cite{Fag} prior to publication and for pointing out ref. \cite{AP}, which simplifies our previous proof of the connectivity lemma.

\end{document}